\documentclass[9pt,twocolumn]{IEEEtran}
\usepackage{amsmath,amssymb,euscript,yfonts,psfrag,latexsym,dsfont,graphicx}
\usepackage{bbm,color,amstext,wasysym,subfig,flushend,parskip,balance}
\graphicspath{{./},{./figures/}}

\newtheorem{thm}{Theorem}

\newtheorem{prop}{Proposition}
\newtheorem{problem}[thm]{Problem}
\newtheorem{remark}[thm]{Remark}

\newcommand{\mR}{{\mathbb R}}

\newcommand{\E}{{\mathbb E}}  

\newcommand{\cR}{{\mathcal R}}
\newcommand{\cS}{{\mathcal S}}

\newcommand{\tr}{\operatorname{trace}}

\newcommand{\trace}{\operatorname{trace}}

\definecolor{grey}{rgb}{0.6,0.6,0.6}
\definecolor{lightgray}{rgb}{0.97,.99,0.99}

\IEEEoverridecommandlockouts
\def\spacingset#1{\def\baselinestretch{#1}\small\normalsize}
\spacingset{.95}

\begin{document}
\title{Optimal control of the state statistics\\ for a linear stochastic system
}

\author{Yongxin Chen, Tryphon Georgiou and Michele Pavon
\thanks{Y.\ Chen and T.T.\ Georgiou are with the Department of Electrical and Computer Engineering,
University of Minnesota, Minneapolis, Minnesota MN 55455, USA; {email: \{chen2468,tryphon\}@umn.edu}}
\thanks{M.\ Pavon is with the Dipartimento di Matematica,
Universit\`a di Padova, via Trieste 63, 35121 Padova, Italy; {email: pavon@math.unipd.it}}
\thanks{The research was supported in part by the AFOSR under Grants FA9550-12-1-0319 and FA9550-15-1-0045}}
\markboth{\today}{}

\maketitle

\begin{abstract}
We consider a variant of the classical linear quadratic Gaussian regulator (LQG) in which penalties on the endpoint state are replaced by the specification of the terminal state distribution. The resulting theory considerably differs from LQG as well as from formulations that bound the probability of violating state constraints.
We develop results for optimal state-feedback control in the two cases where i) steering of the state distribution is to take place over a finite window of time with minimum energy, and ii) the goal is to maintain the state at a stationary distribution over an infinite horizon with minimum power. For both problems the distribution of noise and state are Gaussian. In the first case, we show that provided the system is controllable, the state can be steered to any terminal Gaussian distribution over any specified finite time-interval. In the second case, we characterize explicitly the covariance of admissible stationary state distributions that can be maintained with constant state-feedback control. The conditions for optimality are expressed in terms of a system of dynamically coupled Riccati equations in the finite horizon case and in terms of algebraic conditions for the stationary case. In the case where the noise and control share identical input channels, the Riccati equations for finite-horizon steering become homogeneous and can be solved in closed form. The present paper is largely based on our recent work in \cite{chen2014optimal,CheGeoPav14b} and presents an overview of certain key results.
\end{abstract}


\section{Prologue}
Consider the controlled evolution
\begin{align}\nonumber
dx^u(t)&=A(t)x^u(t)dt+B(t)u(t)dt+B_1(t)dw(t),\\ &\quad x^u(0)=x_0\mbox{ a.s.}\label{controlled}
\end{align}
where $x_0$ is an $n$-dimensional Gaussian vector with density
\begin{equation}\label{initial}\rho_0(x)=(2\pi)^{-n/2}\det (\Sigma_0)^{-1/2}\exp\left(-\frac{1}{2}x'\Sigma_0^{-1}x\right)
\end{equation}
and independent of the standard $p$-dimensional Wiener process $\{w(t)\mid 0\le t\le T\}$; $A$, $B$ and $B_1$ are continuous matrix functions of $t$ taking values in $\mR^{n\times n}$, $\mR^{n\times m}$ and $\mR^{n\times p}$, respectively, $\Sigma_0$ is a symmetric positive definite matrix, and $T<\infty$ represents the end point of a time interval of interest.

The basic paradigm of Linear Quadratic Regulator Theory is to specify a target value for the state vector, e.g., the origin for simplicity, and impose the quadratic penalty
\[
J_T(x(T)):=\E\{x(T)'Mx(T)\}
\]
on possible deviation. Thence, the sum of the quadratic norm of the input
\[
J(u):=\E\left\{\int_0^Tu(t)' u(t) \,dt\right\}
\]
and the terminal cost is sought to be minimized over possible control laws. The solution is well known, see e.g.\ \cite{FR}, \cite[Section II]{separation}. It can obtained by solving the
matrix Riccati equation
\begin{equation}\label{R1}
\dot{\Pi}(t)=-A(t)'\Pi(t)-\Pi(t)A(t)+\Pi(t)B(t)B(t)'\Pi(t)
\end{equation}
for $\{\Pi(t) \mid 0\le t\le T\}$ with boundary condition $\Pi(T)=M$. Here, $M$ is a symmetric matrix which is typically taken to be positive semi-definite, but {\em it does not need to}. 
The variational analysis shows that the optimal feedback control law is
\begin{equation}\label{optcontr}
u(x,t):=-B(t)'\Pi(t)x
\end{equation}
and the optimal cost is
\begin{align*}
J&=J(u^*)+J_T(x^*(T))\\
&= \trace(\Sigma_0\Pi(0))+\int_0^T\trace(B_1(t)B_1(t)'\Pi(t))dt.
\end{align*}

The terminal state $x^*(T)$ of the controlled process
\begin{align}\nonumber
dx^*(t)&=\left(A(t)-B(t)B(t)'\Pi(t)\right)x^*(t)dt+B_1(t)dw(t),\\
&\quad x^*(0)=x_0\mbox{ a.s. }\label{optevolution}
\end{align}
has probability density
 \begin{equation}\label{final0}\rho_T(x)=(2\pi)^{-n/2}\det (\Sigma_T)^{-1/2}\exp\left(-\frac{1}{2}x'\Sigma_T^{-1}x\right),
\end{equation}
where the terminal covariance $\Sigma_T=\Sigma(T)$ is obtained from
\begin{align}\nonumber\dot{\Sigma}(t)&= \left(A(t)-B(t)B(t)'\Pi(t)\right)\Sigma(t)\\&\hspace*{-5pt}+\Sigma(t)\left(A(t)-B(t)B(t)'\Pi(t)\right)'+B_1(t)B_1(t)'\label{sigmaevol}
\end{align}
with boundary condition $\Sigma(0)=\Sigma_0$.
The standard formalism of LQG does not directly address terminal conditions which are only indirectly impacted by the choice of the ``weight" $M$.
In the present work, we study the question of {\em what values of the state covariance $\Sigma(T)$ can be obtained through feedback control}. Moreover, for a given covariance, we consider the problem to {\em steer the system to the corresponding state-distribution with minimal (quadratic) effort for the control input}.


Turning to a time-invariant system, that is \eqref{controlled} with $A,B,B_1$ being independent of time, any state-feedback gain $K$ for which $A-BK$ is a Hurwitz matrix corresponds to a state covariance $\Sigma$ for the stationary controlled state process. It is well known that $\Sigma$ satisfies the Lyapunov equation
\[
(A-BK)\Sigma+\Sigma (A-BK)' +B_1B_1'=0.
\]
In this context we are interested in {\em characterizing positive definite matrices that are admissible as stationary state-covariances}, and to {\em determine a choice for the feedback gain $K$} (assuming that there is more than one) {\em that minimizes power at the control input}.

Both of these problems, to steer and possibly maintain the state statistics of a stochastically driven system, represent variants of the classical {\em regulator problem} which is at the heart of many control applications and entails efficient and accurate steering to a target location,  see e.g.\ \cite{FR,MLS,DD,AT}. The paradigm that is put forth in this paper represents a ``relaxed" version of the classical linear quadratic regulator (LQG) in that hard constraints and penalties on the endpoint state, are replaced by ``soft conditioning'' on the state to be distributed according to a prescribed probability density. This formulation appears natural for applications  in quality control, industrial and manufacturing processes. Further applications may be envisaged in the control of aircrafts, UAVs, and autonomous cars, where specifications on the state distribution are also natural.
Further motivation for this work is provided by modern technological advances that allow us to manipulate micro and thermodynamic systems, and to measure physical properties with unprecedented accuracy. Many such advances rely
heavily on our ability to limit state-uncertainty using feedback, e.g., in oscillators coupled to a heat bath or in steering the collective behavior of a swarm of particles experiencing stochastic forcing. Cutting edge examples include
thermally driven atomic force microscopy,
the control of molecular motors, laser driven reactions, the manipulation of macromolecules,
the ``active cooling'' of devices aimed at measuring gravitational waves,
and the focusing of particle beams
(see \cite{toyabe2010nonequilibrium,gannepalli2005thermally,braiman2003control,hayes2001active,ricci2014low,rowan2000gravitational}).

In order to provide insight on our ability to steer and maintain the probability density of the state, it is useful to recall analogous results regarding controllability of a deterministic system. More specifically,
consider
\begin{equation}\label{eq:linearsystem}
\dot x(t)=A(t)x(t)+B(t)u(t), \;t\in[0,\infty)
\end{equation}
with $A(t),B(t)$ as before, $x(t)\in\mR^n$, $u(t)\in\mR^m$,
and the problem to steer \eqref{eq:linearsystem} from the origin to a given point $x(T)=\xi\in\mR^n$. This is of course  possible for any arbitrary $\xi\in\mR^n$ iff the system is {\em controllable}, i.e.,
the controllability Gramian
\[
\int_0^t \Phi(t,\tau)B(\tau)B(\tau)'\Phi(t,\tau)'d\tau
\]
is positive definite for all $t>0$, or in case $A,B$ are constant, the rank of $[B,\,AB,\ldots,\,A^{n-1}B]$ is $n$; these will be standing assumptions throughout. It is well known that the steering can be effected in a variety of ways, including ``minimum-energy'' control, over any pre-specified interval $[0,T]$. On the other hand, for the case where $A,B$ are time invariant, maintaining a fixed value $\xi$ for the state vector in a stable manner {\em is not always possible}. For this to be the case, using feedback and feedforward control, the equation
\begin{equation}\label{eq:steadystate}
0=(A-BK)\xi + Bu
\end{equation}
must have a solution $(u,K)$ for a constant value of the input $u$ and a corresponding value for $K$ so that $A-BK$ is  Hurwitz (i.e., the feedback system be asymptotically stable).
This reduces to the requirement that
\[
0=A\xi+Bv
\]
for some $v$; if there is such a $v$, we can always choose a suitable $K$ so that $A-BK$ is Hurwitz and then, from $v$ and $K$, we can compute the constant value $u$. Conversely, from $u$ and $K$ we can obtain $v=u-K\xi$.

A similar dichotomy between our ability to steer the distribution of the state vector over a finite window in time, and to maintain an admissible stationary distribution over the infinite time-horizon exists. This is explained next. In what follows, we focus on achieving the control objective with minimal energy and power, respectively, and consider these questions in the case of Gaussian distributions for the noise and state vector. It will be shown that the state-covariance can be assigned at the end of an interval through suitable feedback control if and only if the system is controllable. On the other hand, a positive semidefinite matrix is an admissible stationary state-covariance attained through constant feedback if and only if it satisfies a certain Lyapunov-like algebraic equation.

The ability to specify the mean value of the state-vector reduces to the deterministic problem just discussed. More specifically, since $\E\{x(t)\}=:\bar x(t)$ satisfies \eqref{eq:linearsystem}, controllability or the system is necessary and sufficient to specify $\bar x(T)$ at the end of the interval and this is effected by a deterministic mean value for the input process. Likewise, the mean value for a stationary input must satisfy \eqref{eq:steadystate} to attain $\bar x(t)\equiv\xi$ for a stationary state process. Thus, throughout and without loss of generality we assume that all processes have zero-mean and we only focus on our ability to assign the state-covariance in those two instances.

\section{Finite horizon steering}
Let $\Sigma_T$ be  symmetric and  positive definite $n\times n$ matrix .Consider the ``target'' Gaussian end-point distribution
 \begin{equation}\label{final}
 \rho_T(x)=(2\pi)^{-n/2}\det (\Sigma_T)^{-1/2}\exp\left(-\frac{1}{2}x'\Sigma_T^{-1}x\right),
\end{equation}
 for \eqref{controlled}.
We seek to determine whether is possible to steer (\ref{controlled}) from the initial probability density $\rho_0$ to this ``target'' final probability density $\rho_T$ and, if so, to do this optimally by minimizing
\begin{equation}\label{finiteenergy}
J(u):=\E\left\{\int_0^Tu(t)' u(t) \,dt\right\}<\infty,
\end{equation}
over {\em adapted},
\footnote{I.e., $u\in\mathcal U$ is such that  $u(t)$ only depends on $t$ and on $\{x^u(s)\mid  0\le s\le t\}$ for each $t\in [0,T]$}
{\em finite-energy} control functions such that (\ref{controlled}) has a strong solution on $[0,T]$ and $x^u(T)$ is distributed according to \eqref{final}.
To this end, we let  $\mathcal U$ represent the family of all such admissible control laws and consider the following.

\begin{problem}\label{formalization}  Determine $u^*$ that minimizes $J(u)$ over all $u\in\mathcal U$, i.e., over adapted inputs that steer the system from state-covariance $\Sigma_0$ to $\Sigma_T$.

 As it will be apparent in what follows, our results also cover the case where one or both endpoint marginal densities are delta functions. If $\rho_T(x)=\delta(x-x_0)$ , the optimal control becomes unbounded in suitable directions as $t\nearrow T$ to push the whole diffusion into $x_0$ as in the {\em Brownian bridge}.
Though in the limit, the optimal control may fail to have finite-energy (see also \cite{CheGeo14}).
\end{problem}

\subsection{Conditions for optimality}

Clearly, if $\{\Pi(t) \mid 0\le t\le T\}$ is a solution of the matrix Riccati equation
\eqref{R1}, $u(x,t)$ given by \eqref{optcontr}, and it is such that the end value
 $x^*(T)$ of the Gauss-Markov process
in \eqref{optevolution}
has probability density $\rho_T$, then $u(x^*(t),t)$ is indeed the solution $u^*$ to Problem \ref{formalization}. We recast this observation as follows (Proposition \ref{prop:prop1}).

First, let $\Sigma(t)$ be the state covariance of \eqref{optevolution}, i.e., a solution of \eqref{sigmaevol},
which therefore satisfies the two boundary conditions
\begin{equation}\label{BND}
\Sigma(0)=\Sigma_0, \quad \Sigma(T)=\Sigma_T.
\end{equation}
Then, define
$${\rm H}(t):=\Sigma(t)^{-1}-\Pi(t)
$$
(note that $\Sigma(t)$ is positive definite on $[0,T]$ since, by assumption, $\Sigma_0$ is already positive definite.
A direct calculation using (\ref{sigmaevol}) and (\ref{R1}) leads to \eqref{Schr2} below.
We therefore arrive at a {\em nonlinear} system of equations
\begin{subequations}\label{Schr1234}
\begin{eqnarray}\label{Schr1}
\hspace*{-5pt}\dot{\Pi} &=&-A '\Pi -\Pi A +\Pi B B '\Pi \\
\hspace*{-55pt}\dot{\rm H} &=&-A '{\rm H} -{\rm H} A -{\rm H} B B '{\rm H} \label{Schr2}\\
&&\hspace*{1cm}+\left(\Pi +{\rm H} \right)\left(B B '-B_1 B_1 '\right)\left(\Pi +{\rm H} \right)\nonumber\\
\hspace*{-5pt}\Sigma_0^{-1}&=&\Pi(0)+{\rm H}(0)\label{Schr3}\\
\hspace*{-5pt}\Sigma_T^{-1}&=&\Pi(T)+{\rm H}(T).\label{Schr4}
\end{eqnarray}
\end{subequations}
The case $\Pi(t)\equiv 0$ corresponds to the situation where the uncontrolled evolution already satisfies the boundary marginals and, in that case, ${\rm H}(t)^{-1}$ is simply the prior state covariance.
We summarize as follows.
\begin{prop}\label{prop:prop1} {Assume that} $\{(\Pi(t),{\rm H}(t)) \mid 0\le t\le T\}$ satisfy (\ref{Schr1})-(\ref{Schr4}). Then the feedback control law (\ref{optcontr}) is the solution to Problem \ref{formalization} and the corresponding optimal evolution is given by (\ref{optevolution}).
\end{prop}

For the special case where $B(t)\equiv B_1(t)$, that is, the case where the control inputs and the noise enter the system through identical channels, the nonlinear system of equations \eqref{Schr1234} reduces to~\footnote{This corresponds to a so-called Schr\"odinger system \cite{W}, consisting of a forward and a backward Kolmogoroff (partial) differential equation that are coupled through their boundary conditions. Since the distributions are Gaussian, the system entails matrix differential equations.}
    \begin{subequations}\label{eq:Schr1234}
    \begin{eqnarray}\label{eq:Schr1}
    \hspace*{-5pt}\dot{\Pi} &=&-A '\Pi -\Pi A +\Pi B B '\Pi \\
    \hspace*{-55pt}\dot{\rm H} &=&-A '{\rm H} -{\rm H} A -{\rm H} B B '{\rm H} \label{eq:Schr2}\\
    \hspace*{-5pt}\Sigma_0^{-1}&=&\Pi(0)+{\rm H}(0)\label{eq:Schr3}\\
    \hspace*{-5pt}\Sigma_T^{-1}&=&\Pi(T)+{\rm H}(T).\label{eq:Schr4}
    \end{eqnarray}
    \end{subequations}

For this simplified set of equations \eqref{eq:Schr1234} the existence and uniqueness of solutions follows already from works by Fortet \cite{fortet}, Beurling \cite{Beurling}, Jamison \cite{Jamison}, F\"{o}llmer \cite{W}. Indeed, in this case the problem reduces to a special linear-quadratic version of the so-called Schr\"odinger bridge problem. It should be noted that, while Schr\"odinger bridges \cite{W} are constructed for general diffusions, they are restricted in two essential ways. First the diffusion coefficient is nonsingular (corresponding to $B_1$ being square and invertible) and then, they correspond to the case where control enters in precisely the same way as the noise \cite{dai1991stochastic} ($B=B_1$).

Herein, we are interested in the general case where control and noise excitation channels may differ, albeit, we restrict our attention to linear dynamics and Gaussian statistics.
In \cite{chen2014optimal} we considered \eqref{eq:Schr1234} directly and provided a solution in closed-form.
This solution \cite[Proposition 4]{chen2014optimal}, which is not repeated here due to space limitations,
provides a direct proof of existence of solutions for \eqref{eq:Schr1234} without relying on the theory of the Schr\"odinger bridges. In particular, it establishes {\em feasibility} of Problem \ref{formalization}, i.e.\ that $\mathcal U$ is nonempty and that there exists a minimizer for this special case where $B(t)\equiv B_1(t)$. In general, when $B(t)\neq B_1(t)$ (and time-varying) it is not known whether Problem \ref{formalization} has a minimizer. However, when $B$ and $B_1$ do not depend on time, it turns out that the set of admissible controls $\mathcal U$ is not empty. This is shown below and in Section \ref{sec:num} we also provide an approach that allows constructing suboptimal controls incurring cost that is arbitrarily close to $\inf_{u\in{\mathcal U}}J(u)$.

\subsection{Controllability of state statistics}

We now turn to the ``controllability'' question
of whether there exist admissible controls
to steer the controlled evolution \eqref{controlled}
to a target Gaussian distribution $\Sigma_T$ at the end of a
finite interval $[0,\,T]$.
We do so for the case where
$A,B$ and $B_1$, are time-invariant, and to this end, we
search over controls that can be expressed in state-feedback form
\begin{equation}\label{eq:feedback}
u(x,t)=-K(t) x.
\end{equation}
Then, the state covariance
\[\Sigma(t):=\E\{x(t)x(t)'\}
\]
of \eqref{controlled}
satisfies the Lyapunov differential equation
\begin{equation}\label{eq:covariancedynamics}
\dot\Sigma(t)=(A-BK(t))\Sigma(t)+\Sigma(t)(A-BK(t))'+B_1B_1'
\end{equation}
with $\Sigma(0)=\Sigma_0$. Regardless of the choice of
$K(t)$, \eqref{eq:covariancedynamics} specifies dynamics that leave invariant the cone of positive semi-definite symmetric matrices
\[
\cS_n^+:=\{\Sigma \mid \Sigma\in\mR^{n\times n},\;\Sigma=\Sigma'\geq 0\}.
\]
To see this, note that the solution to \eqref{eq:covariancedynamics} is of the form
\[
\Sigma(t)=\hat\Phi(t,0)\Sigma_0\hat\Phi(t,0)' +
\int_0^t\hat\Phi(t,\tau)B_1B_1'\hat\Phi(t,\tau)' d\tau
\]
where $\hat\Phi(t,0)$ satisfies
\[
\frac{\partial \hat\Phi(t,0)}{\partial t}=(A-BK(t))\hat\Phi(t,0)
\]
and $\hat\Phi(0,0)=I$, the identity matrix; i.e., $\hat\Phi(t,0)$ is the state-transition matrix of the system $\dot x(t)=(A-BK(t))x(t)$.

 Our interest is in our ability to specify $\Sigma(T)$ via a suitable choice of $K(t)$.
To this end, define
\[
U(t):=-\Sigma(t)K(t)',
\]
and observe that $U(t)$ and $K(t)$ are in bijective correspondence provided that $\Sigma(t)>0$ (which follows from $\Sigma_0>0$). Thus, we now consider the differential Lyapunov system
\begin{equation}\label{eq:diffLyapunov}
\dot\Sigma(t)=A\Sigma(t)+\Sigma(t)A'+BU(t)'+U(t)B'.
\end{equation}
Reachability/controllability of a differential system such as \eqref{eq:linearsystem}, or \eqref{eq:diffLyapunov}, is the property that with suitable bounded control input $u(t)$, or $U(t)$, respectively, the solution can be driven to any finite value. Interestingly, if either \eqref{eq:linearsystem} or \eqref{eq:diffLyapunov} is controllable, so is the other. But, more importantly, when \eqref{eq:diffLyapunov} is controllable, the control authority allowed is such that steering from one value for the covariance to another can be done by remaining within the non-negative cone. This is stated as next (see \cite[Theorem 3]{CheGeoPav14b}).

\begin{thm}\label{thm:thm1}
The Lyapunov system \eqref{eq:diffLyapunov} is controllable iff $(A,B)$ is a controllable pair. Furthermore, if \eqref{eq:diffLyapunov} is controllable, given any two positive definite matrices $\Sigma_0$ and $\Sigma_T$ and an arbitrary $Q\geq 0$, there is a smooth input $U(\cdot)$ so that the solution of
the (forced) differential equation
\begin{equation}\label{eq:diffLyapunov2}
\dot\Sigma(t)=A\Sigma(t)+\Sigma(t)A'+BU(t)'+U(t)B' +Q
\end{equation}
satisfies the boundary conditions $\Sigma(0)=\Sigma_0$ and $\Sigma(T)=\Sigma_T$
and $\Sigma(t)>0$ for all $t\in[0,T]$.
\end{thm}
\begin{remark}
Although the result on the controllability of state covariance we establish is for time-invariant system, we believe that a similar result holds for time-varying system.
\end{remark}


\subsection{Finite interval minimum energy steering of state statistics}\label{sec:num}
We are interested in computing an optimal choice for feedback gain $K(t)$ so that
the control $u(t)=-K(t)x(t)$ steers \eqref{controlled}
from an initial state-covariance $\Sigma_0$ at $t=0$ to the final $\Sigma_T$ at $t=T$. The expected control-energy functional
    \begin{eqnarray}\label{eq:functional}
       J(u)&:=& \E\left\{\int_0^T u(t)'u(t)dt\right\}\\
       &=&\int_0^T \tr(K(t)\Sigma(t)K(t)')dt\nonumber
    \end{eqnarray}
needs to be optimized over $K(t)$ so that \eqref{eq:covariancedynamics} as well as the boundary conditions
\begin{subequations}\label{eq:sdpproblem}
\begin{equation}\label{eq:boundary}
\Sigma(0)=\Sigma_0, \mbox{ and }\Sigma(T)=\Sigma_T
\end{equation}
hold.

If instead of the choice $K(t)$ we sought to optimize over $U(t):=-\Sigma(t)K(t)'$ and $\Sigma(t)$, the functional \eqref{eq:functional} becomes
\[
J=\int_0^T \tr(U(t)'\Sigma(t)^{-1}U(t))dt
\]
which is jointly convex in $U(t)$ and $\Sigma(t)$, while \eqref{eq:covariancedynamics} is replaced by
\begin{equation}\label{eq:diffeqB1U}
\dot\Sigma(t)=A\Sigma(t)+\Sigma(t) A'+BU(t)'+U(t)B'+B_1B_1'
\end{equation}
which is now linear in both! As a consequence, the optimization can be written in the form of a semi-definite program to minimize
\begin{equation}\label{eq:sdp}
 \int_0^T \tr(Y(t))dt
 \end{equation}
 subject to (\ref{eq:boundary}-\ref{eq:diffeqB1U}) and
 \begin{equation}
\left[\begin{matrix}Y(t)& U(t)' \\U(t) & \Sigma(t)\end{matrix}\right]\ge 0.
\end{equation}
\end{subequations}
This can be solved numerically after discretization in time, and a (suboptimal) gain
recovered using the correspondence $K(t)=-U(t)'\Sigma(t)^{-1}$.

\section{Infinite-horizon steering}

When $A$, $B$ and $B_1$ do not depend on time,
we seek a constant state feedback law $u(t)=-Kx(t)$
to maintain a stationary state-covariance $\Sigma>0$ for \eqref{controlled}. In particular, we are interested
in one that minimizes the expected input power (energy rate)
\begin{eqnarray}\label{eq:power}
J_{\rm power}(u)&:=&\E\{u'u\}
\end{eqnarray}
and thus we are led to the following problem\footnote{ An equivalent problem is to minimize
$\lim_{T\to \infty}\frac{1}{T}\E\left\{\int_0^Tu(t)'u(t)dt\right\}$
for a given terminal state covariance as $T\to\infty$.}.
\begin{problem}\label{problem2} Determine $u^*$ that minimizes $J_{\rm power}(u)$ over all $u(t)=-Kx(t)$ such that
\begin{equation}\label{feedbackdynamics}
dx(t)=(A-BK)x(t)dt+B_1dw(t)
\end{equation}
admits
\begin{equation}\label{invdensity}
\rho(x)=(2\pi)^{-n/2}\det (\Sigma)^{-1/2}\exp\left(-\frac{1}{2}x'\Sigma^{-1}x\right)
\end{equation}
as invariant probability density.
\end{problem}
In general, this problem may not have a solution. In Theorem \ref{admissiblestate3} we provide conditions that ensure $\Sigma$ is admissible as a stationary state covariance for a suitable input.
 Moreover, as it will be apparent from what follows, even when the problem is feasible, i.e., there exist controls which maintain $\Sigma$, an optimal control may fail to exist. The problem has connections to Jan Willems' classical work on the Algebraic Riccati Equation \cite{willems1971least} and this is discussed following Proposition \ref{prop:prop2}.

\subsection{Condition for optimality}

Let us start by observing that the problem admits the following finite-dimensional reformulation. Let $\mathcal K$ be the set of all $m\times n$ matrices $K$ such that the corresponding feedback matrix $A-BK$ is Hurwitz. Observe that
\[
\E\{u'u\}=\E\{x'K'Kx\}=\tr(K\Sigma K')
\]
Then Problem \ref{problem2} reduces to finding a $m\times n$ matrix $K\in\mathcal K$ which minimizes
\begin{equation}\label{criterion}
J(K)=\tr\left(K\Sigma K'\right)
\end{equation}
subject to the constraint
\begin{equation}
(A-BK)\Sigma+\Sigma(A'-K'B')+B_1B_1'=0.\label{constraint}
\end{equation}
Consider the Lagrangian
 \begin{eqnarray}
 \mathcal{L}(K,\Pi)&=&\tr\left(K\Sigma K'\right)\\\nonumber&&\hspace*{-1.5cm}+\tr\left(\Pi((A-BK)\Sigma+\Sigma(A'-K'B')+B_1B_1')\right).
 \end{eqnarray}
Note that since $\mathcal K$ is {\em open}, a minimum point may fail to exist.
Standard variational analysis \cite{CheGeoPav14b} leads to the form
$
K=B'\Pi
$
for the optimal gain. We summarize this as follows.

\begin{prop}\label{prop:prop2}
Assume that there exists a symmetric matrix $\Pi$ such that $A-BB'\Pi$ is a Hurwitz matrix and
\begin{equation}\label{sigmastat}
(A-BB'\Pi)\Sigma+\Sigma(A-BB'\Pi)'+B_1B'_1=0
\end{equation}
holds.
Then, the solution to Problem \ref{problem2} is
\begin{equation}\label{statoptcontr}
u^*(t)=-B'\Pi x(t)
\end{equation}
\end{prop}


We now draw a connection to some classical results due to  Jan Willems \cite{willems1971least}.  In our setting, minimizing \eqref{eq:power} is equivalent to minimizing
\begin{equation}\label{eq:alternate}
J_{\rm power}(u)+\E\{x'Qx\}
\end{equation}
for an arbitrary symmetric matrix $Q$ since the portion
\[
\E\{x'Qx\}=\tr\{Q\Sigma\}
\]
is independent of the choice of $K$ (because $\Sigma$ is given).
On the other hand, minimization of \eqref{eq:alternate} for specific $Q$, but without the constraint that $\E\{xx'\}=\Sigma$, was studied by Willems \cite{willems1971least} and is intimately related to the {\em maximal} solution of the Algebraic Riccati Equation (ARE)
\begin{equation}\label{ARE}
A'\Pi+\Pi A-\Pi BB'\Pi+Q=0.
\end{equation}
Specifically, under the assumption that the Hamiltonian matrix
$$H=\left[\begin{matrix} A &-BB'\\-Q & -A'\end{matrix}\right]
$$
has no pure imaginary eigenvalues, Willems' result states  that $A-BB'\Pi$ is Hurwitz and that \eqref{statoptcontr} is the optimal solution.
Thus, starting from a symmetric matrix $\Pi$ as in Proposition \ref{prop:prop2}, we can define $Q$ using
\[
Q=-A'\Pi-\Pi A+\Pi BB'\Pi.
\]
By Willems' results, \eqref{ARE} has at most one ``stabilizing'' solution $\Pi$, and therefore, the matrix in the proposition coincides with the maximal solution to \eqref{ARE}. Therefore, if Problem~\ref{problem2} has a solution, this corresponds to the maximal solution of the ARE~\eqref{ARE} for a particular choice of $Q$.
Interestingly, neither $\Pi$ nor $Q$, corresponding to an optimal control law for which \eqref{sigmastat} holds, are unique, whereas $K$ is. The computation and the uniqueness of the optimal gain $K$ will be discussed later on in Section \ref{sec:mestationary}.

\subsection{Characterization of stationary state statistics}

We now focus on what values of the state covariance are admissible in that they can be obtained by state feedback.

The constraint \eqref{constraint} already implies that
\begin{subequations}\label{eq:equivalent}
\begin{eqnarray}\label{eq:lyapunov3}
&&A\Sigma + \Sigma A'+B_1B_1'+BX'+XB'=0\\
&&\mbox{can be solved for $X$}.\nonumber
\end{eqnarray}
In particular we can take $X=-\Sigma K'$ in \eqref{constraint}. Thus, the solvability of \eqref{eq:lyapunov3} is a necessary condition for $\Sigma$ to qualify as a stationary state-covariance attained via feedback.
This condition, \eqref{eq:lyapunov3}, can be written as a rank condition:
\begin{equation}\label{eq:rank2}
{\rm rank}\left[\begin{matrix}
A\Sigma+\Sigma A'+B_1B_1' & B\\
B & 0
\end{matrix}\right]
=
{\rm rank}\left[\begin{matrix}
0 & B\\
B & 0
\end{matrix}\right],
\end{equation}
which ensures that $A\Sigma + \Sigma A'+B_1B_1'$ is in the range of the linear map $X\mapsto BX'+XB'$, cf.\  \cite[Proposition 1]{georgiou2002structure}.
\end{subequations}

Conversely, given $\Sigma>0$ that satisfies \eqref{eq:equivalent} and a solution $X$ of \eqref{eq:lyapunov3},
then \eqref{constraint} holds with $K=-X'\Sigma^{-1}$. Therefore, provided $A-BK$ is a Hurwitz matrix, $\Sigma$ {\em is an admissible stationary covariance}. Now, the property of $A-BK$ being Hurwitz can be guaranteed when
$(A-BK,\,B_1)$ is a controllable pair. In turn, controllability of $(A-BK,\,B_1)$ is guaranteed when $\cR(B)\subseteq \cR(B_1)$. Thus, we have established the following.

\begin{thm}\label{admissiblestate3}
Consider the Gauss-Markov model \eqref{controlled}
and assume that $\cR(B)\subseteq \cR(B_1)$. A positive-definite matrix $\Sigma$ can be assigned as the stationary state covariance via a suitable choice of state-feedback if and only if $\Sigma$ satisfies any of the equivalent statements (\ref{eq:lyapunov3}-\ref{eq:rank2}). \end{thm}

Holtz and Skelton  \cite{hotz1987covariance} considered the problem to maintain a stationary state covariance with feedback and provide as condition that $A\Sigma + \Sigma A'+B_1B_1'$ is in the null space of the operator
$
Y\mapsto \Pi_{{\mathcal R}(B)^\perp}Y\Pi_{{\mathcal R}(B)^\perp}$. It turns out that this condition is equivalent to (\ref{eq:lyapunov3}-\ref{eq:rank2}), see \cite{CheGeoPav14b}.
On the other hand, conditions (\ref{eq:lyapunov3}-\ref{eq:rank2}) were obtained in \cite{georgiou2002structure,georgiou2002spectral}, for
the special case when $B=B_1$, as being necessary and sufficient for a positive-definite matrix to materialize as the state covariance of the system driven by a stationary stochastic process (not-necessarily white).

\subsection{Minimum energy control}\label{sec:mestationary}

As we just noted, a positive definite matrix $\Sigma$ is admissible as a stationary state-covariance provided \eqref{eq:lyapunov3} holds for some $X$ and $A+BX'\Sigma^{-1}$ is a Hurwitz matrix. The condition $\cR(B)\subseteq \cR(B_1)$ is a sufficient condition for the latter to be true always, but it may be true even if $\cR(B)\subseteq \cR(B_1)$ fails (see Example 1 in Section \ref{sec:example}). Either way, the expected input power (energy rate)
\begin{eqnarray}
\E\{u'u\}
&=&\tr(K\Sigma K')
\\\nonumber
&=&\tr(X'\Sigma^{-1}X)
\end{eqnarray}
in either $K$, or $X$.
Thus, assuming that $\cR(B)\subseteq \cR(B_1)$ holds, and in case \eqref{eq:lyapunov3} has multiple solutions, the optimal constant feedback gain $K$ can be obtained by solving the convex optimization problem
\begin{equation}\label{eq:XSX}
\min\left\{\tr(X'\Sigma^{-1}X)\mid \mbox{ \eqref{eq:lyapunov3} holds } \right\}.
\end{equation}
\begin{remark}
In case $\cR(B)\not\subseteq \cR(B_1)$, the condition that $A-BK$ be Hurwitz needs to be verified separately. If this fails, we cannot guarantee that $\Sigma$ is an admissible stationary state-covariance that can be maintained with constant state-feedback. However, it is always possible to maintain a state-covariance that is arbitrarily close. To see this, consider the control
\[
K_\epsilon = K + \frac12 \epsilon B' \Sigma^{-1}
\]
for $\epsilon>0$. Then,
\begin{eqnarray*}
(A-BK_\epsilon)\Sigma + \Sigma(A-BK_\epsilon)'&=&-\epsilon BB'-B_1B_1'\\
& \leq& -\epsilon BB'.
\end{eqnarray*}
The fact that $A-BK_\epsilon$ is Hurwitz is obvious. If now $\Sigma_\epsilon$ is the solution to
\begin{eqnarray*}
(A-BK_\epsilon)\Sigma_\epsilon + \Sigma_\epsilon (A-BK_\epsilon)'&=&-B_1B_1'
\end{eqnarray*}
the difference $\Delta=\Sigma-\Sigma_\epsilon\geq 0$ and satisfies
\begin{eqnarray*}
(A-BK_\epsilon)\Delta + \Delta(A-BK_\epsilon)'&=&-\epsilon BB',
\end{eqnarray*}
and hence is of $o(\epsilon)$.
\end{remark}

\section{Examples}\label{sec:example}
\subsection*{Example 1}
Consider particles that are modeled by
   \begin{eqnarray*}
       dx(t) &=& v(t)dt + dw(t)\\
       dv(t) &=& u(t)dt.
   \end{eqnarray*}
Here, $u(t)$ is the control input (force) at our disposal, $x(t)$ represents position and $v(t)$ velocity (integral of acceleration due to input forcing), while $w(t)$ represents random displacement due to impulsive accelerations. Alternatively, $\int^t v(\tau)d\tau$ may represent actual position while $x(t)$ noise measurement of position.

The purpose of the example is to highlight a case where the control is handicapped compared to the effect of noise. Indeed, the displacement $w(t)$ directly affects $x(t)$ while the control effort $u(t)$ needs to be integrated before it mitigates the effect of $w(t)$ on the position $x(t)$ of the particles.

Another interesting fact that this example highlights is that $\cR(B)\not\subseteq \cR(B_1)$ is not necessary
for being able to ensure that, particular positive definite $\Sigma$'s in the range of $X\mapsto BX'+XB'$,  can be maintained as stationary state covariances.
Indeed, here, $\cR(B)\not\subseteq \cR(B_1)$ since
$B=[0,~1]'$ while $B_1=[1,~0]'$. Yet, if we choose
\begin{equation}\label{eq:stationaryvalue}
\Sigma_1=\left[\begin{matrix}1&-1/2\\ \hspace*{2pt}-1/2& \phantom{-}1/2\end{matrix}\right]
\end{equation}
as a candidate stationary state-covariance, it can be seen that
\eqref{eq:lyapunov3} has a unique solution $X$ giving rise to  $K=\left[1,~1\right]$ and a stable feedback since $A-BK$ is Hurwitz.

\begin{figure}\begin{center}
\includegraphics[width=0.47\textwidth]{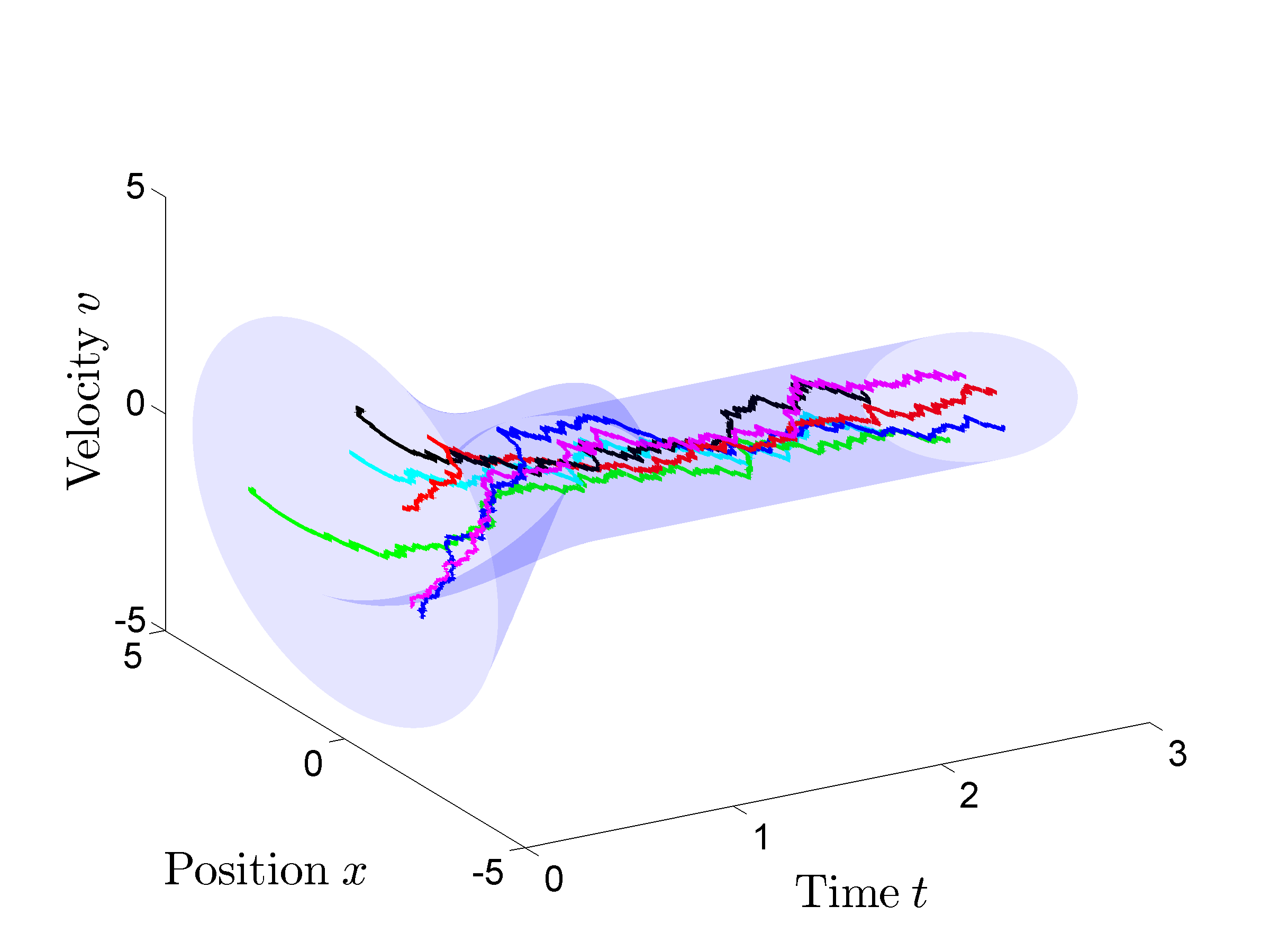}
   \caption{State trajectories (Example 1)}
   \label{fig:Eg1Phase1}
\end{center}\end{figure}

\begin{figure}\begin{center}
\includegraphics[width=0.47\textwidth]{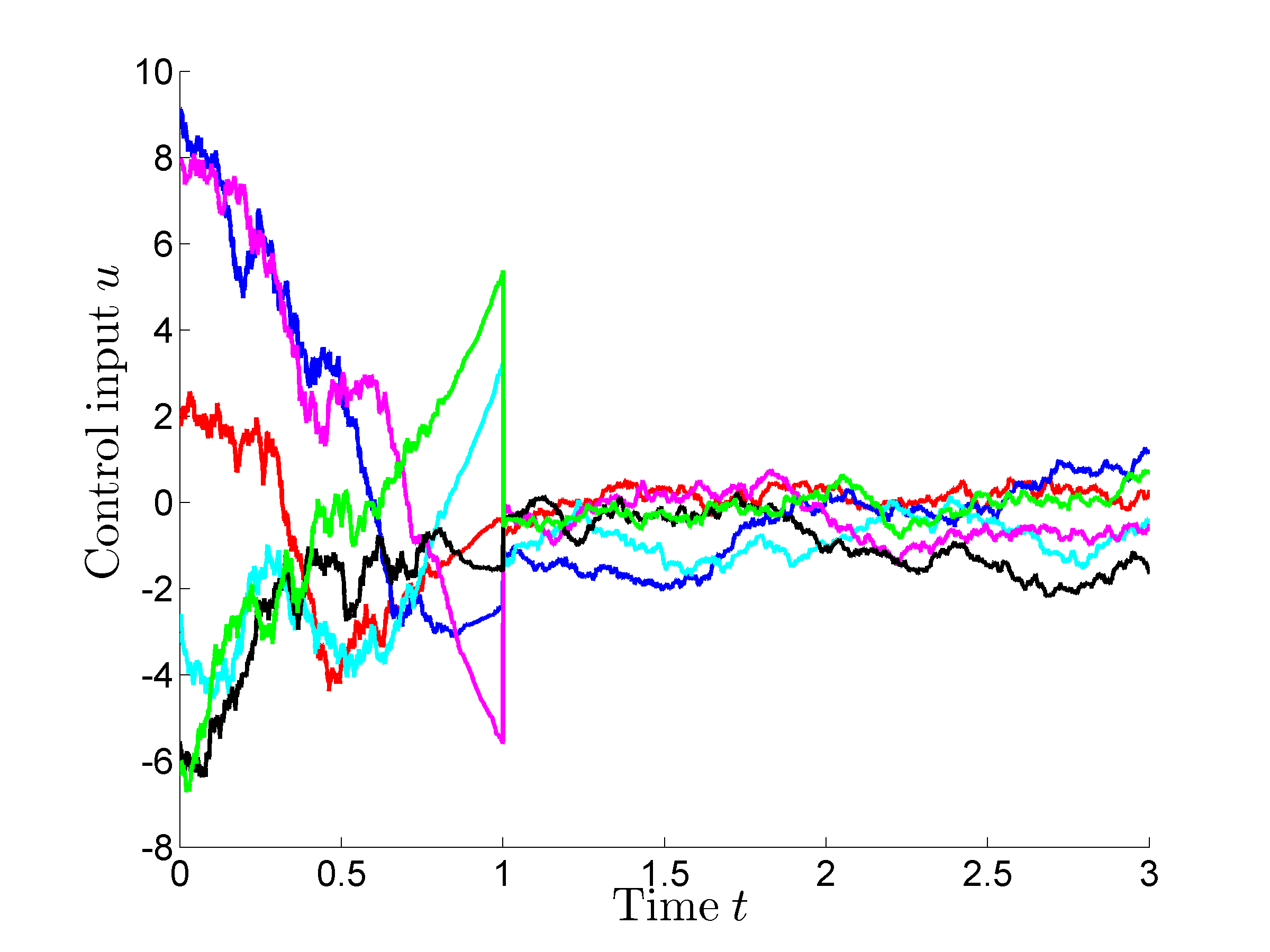}
   \caption{Control inputs (Example 1)}
   \label{fig:Eg2Control1}
\end{center}\end{figure}

We now wish to steer the spread of the particles from an initial Gaussian distribution with $\Sigma_0=2I$ at $t=0$ to the terminal marginal $\Sigma_1$
at $t=1$, and from there on, since $\Sigma_1$ is an admissible stationary state-covariance, to maintain with constant state-feedback control. As we just noted, this is indeed possible.
Figure~\ref{fig:Eg1Phase1} depicts typical sample paths in phase space, as a function of time, that are attained using i) the suboptimal feedback strategy derived following \eqref{eq:sdpproblem} over the time interval $[0,\,1]$ and, ii) static state-feedback with $K=[1,\,1]$ over the time window $[1,\,2]$. Figure \ref{fig:Eg2Control1} displays the corresponding control action for each trajectory over the complete time interval $[0,\,2]$, which consists of the ``transient'' interval $[0,\,1]$ to the target distribution and the ``stationary'' interval $[1,\,2]$.

\subsection*{Example 2}
We now consider a more detailed model of inertial particles where random acceleration as well as control are effected through forcing along the same channel. Yet, the control is again ``handicapped'' due to lag in actuation dynamics. An additional aspect of the example is that it highlights a case where ``target state statistics'' may be prescribed for only a portion of the combined controller-system dynamics, namely, for the dynamics corresponding to the inertial particles and not the controller.

To this end, we consider the second-order process
\begin{eqnarray}\nonumber
dx(t) &=& \phantom{-}v(t)dt\\\nonumber
dv(t) &=& \phantom{-}x_c(t)dt + dw(t)\\
dx_c(t) &=&-x_c(t)dt + u(t)dt.\label{eq:third order}
\end{eqnarray}
Here, $x_c$ is the 1-dimensional state/output of the actuator and represents force, $u(t)$ represents the control signal to the actuator while, once again, $w(t)$ represents random displacement due to impulsive accelerations.
We (arbitrarily) select
\[
\Sigma_{x,v}=\left[\begin{matrix}7/4&0\\ \hspace*{2pt}0&3/4\end{matrix}\right]
\]
as a desirable steady-state covariance for the projection of the process onto the phase-plane of the particle dynamics $(x,v)$. (The additional component $x_c$ corresponding to the actuator is not shown.) First, we need to determine whether $\Sigma_{x,v}$ is indeed an admissible steady-state covariance and, if so, to determine an optimal choice for the constant state-feedback gain that ensures the state is distributed accordingly. To this end, we seek a choice of a variance $\Sigma_{x_c}$ for $x_c$ and of a cross-covariance $Y$ between $x_c$ and $(x,v)$ so that
\[
\Sigma=\left[\begin{matrix}\Sigma_{x,v} & Y\\Y' & \Sigma_{x_c}\end{matrix}\right]>0
\]
is a admissible stationary state-covariance for \eqref{eq:third order}. For this to be true, we need to verify that \eqref{eq:rank2} holds with
\[A=\left[\begin{matrix}0&1&0\\0&0&1\\0&0&-1\end{matrix}\right],\;
B=\left[\begin{matrix}0\\0\\1\end{matrix}\right],\,B_1=\left[\begin{matrix}0\\1\\0\end{matrix}\right].
\]
This is indeed the case for $Y=[-3/4 ~-1/2]', ~\Sigma_{x_c}=3/4$.

\begin{figure}\begin{center}
\includegraphics[width=0.47\textwidth]{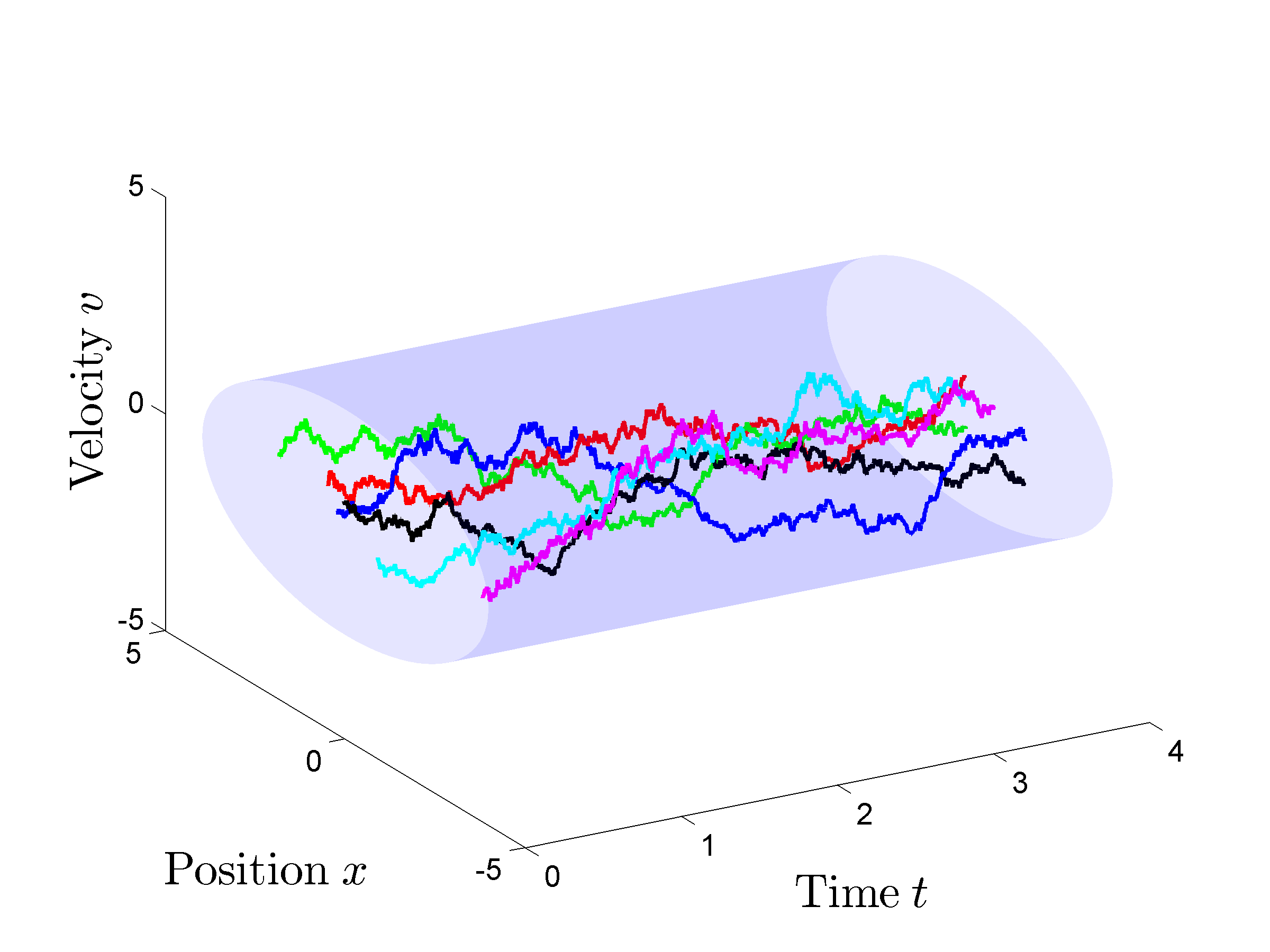}
   \caption{Steady state trajectories in phase space (Example 2)}
   \label{fig:Eg3phase1}
\end{center}\end{figure}

For the above choice of $\Sigma$,
the optimal gain and power are found to be $K=[1 ~3~2]$ and $J=5/2$ by solving \eqref{eq:XSX} using e.g., \cite{cvx}. For this solution, stationary state trajectories, projected onto the $(x,v)$-coordinates are now displayed in Figure \ref{fig:Eg3phase1}.

\begin{figure}\begin{center}
\includegraphics[width=0.47\textwidth]{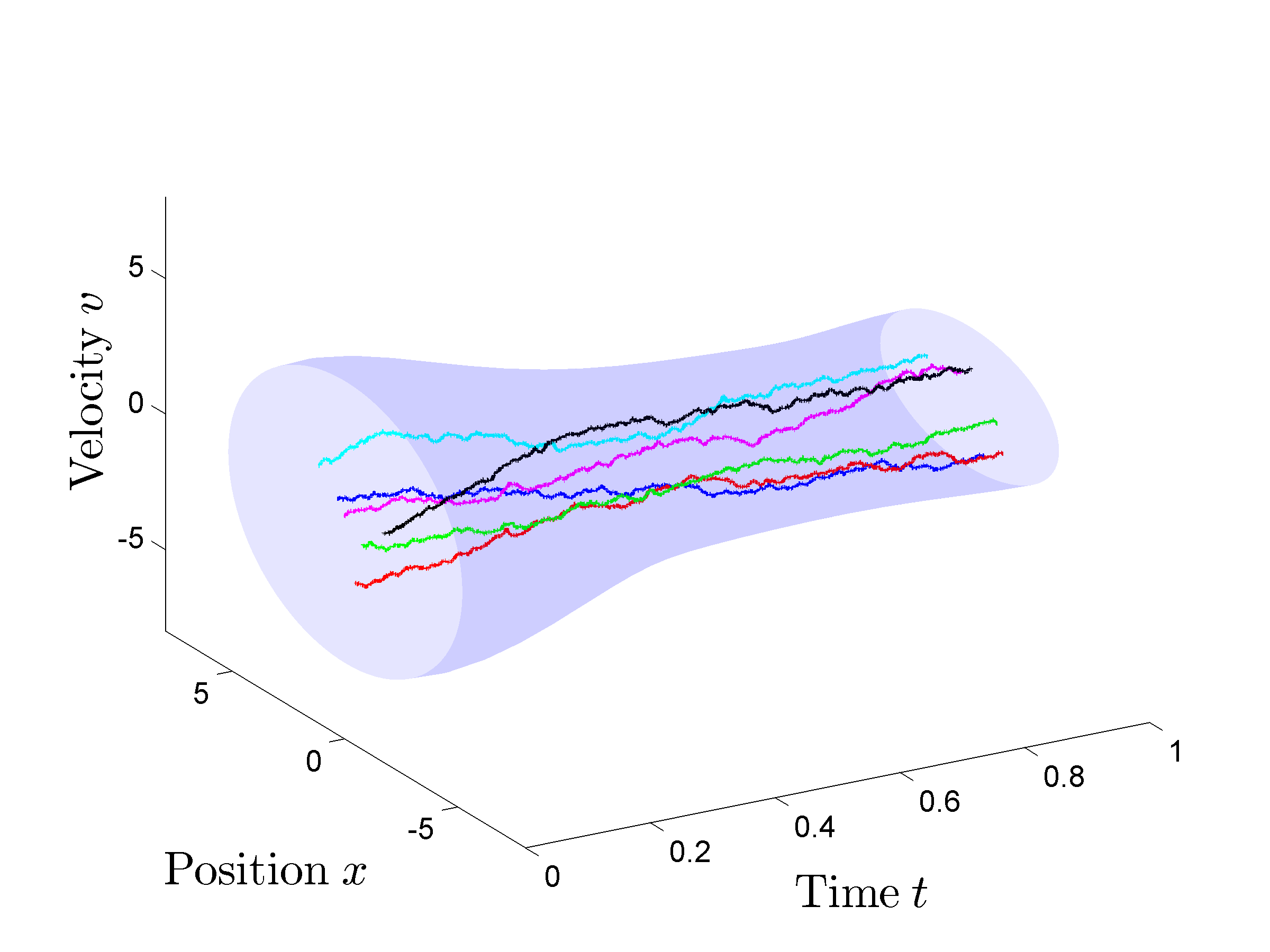}
   \caption{Finite-interval steering in phase space (Example 2)}
   \label{fig:Eg3phase2}
\end{center}\end{figure}

As before the steering between specified Gaussian probability densities
over an interval $[0,1]$ follows Section \ref{sec:num}. For completeness,
we display in Figure \ref{fig:Eg3phase2} sample paths corresponding to the
transition between marginals with covariance matrices
$\Sigma_0=3I$ to $\Sigma_1=\Sigma$, respectively. The figure shows the projection onto
the $(x,v)$-component of the process that corresponds to position and velocity.

\section{Epilogue}
In this paper, we have outlined some of the main results in \cite{chen2014optimal,CheGeoPav14b}. These results, in the finite horizon problem, may be viewed as a generalisation of the theory of Schr\"{o}dinger bridges \cite{W} for Gauss-Markov models to the situation where the noise does not affect all directions and where the control also acts through a different channel (recent attempts to address linear stochastic systems were also limited to non-degenerate diffusions where the control and noise channel are identical \cite{beghi1997continuous,vladimirov2012minimum}). Although the theory is far from complete, our work has provided the first {\em implementable form} of the optimal control. Taking into consideration also the algorithmic solutions of the Schr\"{o}dinger system provided in \cite{GP}, we have now at disposal a new powerful tool to attack a variety of problems in classical and modern control  including some applications in unexpected directions such as the {\em optimal mass transport problem} \cite{CheGeoPav14c,CheGeoPav15a}.

{
\bibliographystyle{IEEEtran}
\bibliography{./refs}
}
\end{document}